\documentclass[11pt,reqno]{article}

\title {Transience of percolation clusters on wedges}

\author{Omer Angel \and Itai Benjamini \and Noam Berger \and Yuval Peres}

\usepackage{latexsym}
\usepackage{amsfonts}
\usepackage{amsmath}
\usepackage{amsthm}

\newtheorem{claim}{Claim}
\newtheorem{thm}{Theorem}
\newtheorem{defn}[thm]{Definition}
\newtheorem{lemma}{Lemma}
\newtheorem{cor}[thm]{Corollary}

\newtheorem{remark}{Remark}

\newcommand{\R}{\mbox{$\bf R$}}
\newcommand{\E}{\mbox{$\bf E$}}
\newcommand{\ep}{\varepsilon}

\newcommand{\N}{{\mathbb N}}
\newcommand{\Z}{{\mathbb Z}}
\newcommand{\W}{{\mathcal W}}
\newcommand{\V}{{\mathcal V}}
\newcommand{\G}{{\mathcal G}}
\renewcommand{\P}{{\bf P}}

\newcommand{\len}{{\text{Length}}}
\newcommand{\olle}{{H\"{a}ggstr\"{o}m }}
\newcommand{\ve}{\text{ and }}
\newcommand{\diam}{\text{diam}}

\begin{document}

\maketitle

\begin{abstract}
We study random walks on supercritical percolation
clusters on wedges in $\Z^3$, and show that the infinite percolation
cluster is (a.s.) transient whenever the wedge is transient.
This solves a question raised by O. \olle and E. Mossel.
We also show that for convex gauge functions satisfying
a mild regularity condition, the existence of a finite energy flow on $Z^2$
is equivalent to the (a.s.) existence of a finite energy flow on the
supercritical  percolation cluster. This solves a question of C. Hoffman.
\end{abstract}



\section{Introduction}
For simple random walk in the  $\Z^d$ lattice,
Polya showed in 1920 that the transition from recurrence
to transience occurs when $d$ increases from $2$ to $3$.
The transition boundary is more sharply delineated
by a 1983 result of T. Lyons, concerning {\it wedges}.
For an
increasing positive function $h$, the wedge $\W_h$ is the
subgraph of $\Z^3$ induced by the vertices
\begin{equation*}
V(\W_h) = \left\{(x,y,z) \mid x\geq0 \text{ and } |z|\leq h(x)\right\} \,.
\end{equation*}
T. Lyons \cite{Ly} proved that $\W_h$ is transient if and only if
\begin{equation} \label{eq:lyons}
\sum_{j=1}^\infty  \frac{1}{jh(j)}  <  \infty.
\end{equation}
(A locally finite graph is called transient or recurrent according to
the type of simple random walk on it.)
It is well-known (see, e.g., \cite{DS} or \cite{climb}) that if $G$ is
recurrent then so is any
subgraph of $G$.

In \cite{GKZ} Grimmett, Kesten and Zhang proved that the infinite cluster
of supercritical percolation in $\Z^d$ is transient when $d\geq 3$.
A different proof and some extensions were given in \cite{BPP}.
\olle and Mossel~\cite{EO} sharpened the methods of \cite{BPP} and
 showed that if the increasing positive function $h$ satisfies
\begin{equation}\label{hagmos}
\sum_{j=1}^\infty  \frac{1}{j \sqrt{h(j)}}  <  \infty \,,
\end{equation}
 then the infinite cluster
of supercritical percolation in $\W_h$ is transient.

The condition (\ref{hagmos}) is strictly stronger than
 Lyons' condition (\ref{eq:lyons});
In particular,  the function $h(j)=h_r(j)=\log^{r}(j)$ satisfies
(\ref{hagmos})  iff $r>2$, while it satisfies
(\ref{eq:lyons})  for all $r>1$.
\olle and Mossel asked what is the type of percolation clusters in
wedges  $\W_h$ where $h$ satisfies Lyons' condition
(\ref{eq:lyons}) , but does not satisfy
condition (\ref{hagmos}).

Our main result answers this question:
\begin{thm} \label{thm:main}
Let $h$ be a positive increasing function.
The infinite cluster of supercritical percolation on the wedge $\W_h$ is
transient if and only if $\W_h$ is transient,
i.e., if and only if $h$ satisfies (\ref{eq:lyons}).
\end{thm}

\medskip

A useful notion will be that of the {\it core\/} of a subgraph.
The core consists of those vertices which are in some sense
far from the boundary:

\begin{defn}
Let $A$ be a subgraph of a graph $G$, and $v_0$ be a fixed vertex in $G$.
{}For $C>0$ define the $C$-core of $A$ to be the subgraph of $A$ induced by
\begin{equation*}
\left\{ v\in A | d(v,\partial A) > C \log d(v,v_0) \right\},
\end{equation*}
where $d$ denotes the graph metric on $G$.
\end{defn}

Theorem \ref{thm:main} is a consequence of the following more general
statement:

\begin{thm} \label{thm:core}
Let $d\geq3$ and let $p>p_c(\Z^d)$. There exists a constant
$C=C_{d,p}$ with the following property: Let $A$ be a subgraph of
$\Z^d$ s.t.
  (I) the $C$-core of $A$ is transient, and
  (II) $p$-percolation on $A$ has, a.s., a unique infinite cluster.
Then, a.s., the infinite $p$-percolation cluster of $A$ is transient.
\end{thm}

From this we infer the following:
\begin{cor}
Let $A$ be a subgraph of $\Z^d$ such that for any $C$ its $C$-core is
transient and such that any supercritical percolation on $A$ has a.s. a
unique infinite cluster, then this infinite cluster is a.s. transient.
\end{cor}

We will restate and prove Theorem~\ref{thm:core} as
a theorem about flows. Consider each
undirected edge of a graph $G$ as two directed edges, one in each
direction. Let $vw$ be the directed edge from $v$ to $w$. A {\bf flow} $F$
on $G$ with source $v_0$ is an edge function such that $F(vw)=-F(wv)$ and
such that for any vertex $v \not= v_0$:
$\, \sum_{w} F(vw) = 0 $. Recall the following definitions (\cite{Hoff}).

\begin{defn}
Let $g:\R\to\R$ be a function. The {\bf $g$-energy} of a flow $F$ on a graph,
denoted by $H_g(F)$ is defined to be
\begin{equation*}
\sum_{e\in E}{g(|F(e)|)}
\end{equation*}
\end{defn}

{}For $d$ and $\alpha$ we define $\Psi_{d,\alpha}$ to be
\begin{equation*}
\Psi_{d,\alpha}(x)= \frac {|x|^{d/(d-1)}} {\log (1+|x|^{-1})^\alpha}
\end{equation*}.

\begin{defn}
The {\bf $(d,\alpha)$-energy} of a flow $F$ on a graph, denoted
$H_{d,\alpha}(f)$ is the $\Psi_{d,\alpha}$-energy of $F$, i.e.
\begin{equation*}
H_{d,\alpha}(f) = \sum_{e\in E} \Psi_{d,\alpha}(F(e))
\end{equation*}
\end{defn}

It is well known that a
 graph is transient if and only if it has a flow with finite
$(2,0)$-energy; see  (see \cite{Ly}, \cite{DS}, or \cite{climb}).
It is also known that on $\Z^d$ there are flows with finite $(d,1+\ep)$-energy
but no flow with finite $(d,1)$-energy (see \cite{Ly} and \cite{HoMo}).
In \cite{HoMo} Hoffman and Mossel
(refining an earlier result of Levin and Peres~\cite{LP})
proved that the same is
true for the infinite cluster of super critical percolation in $\Z^d$,
provided that $d \geq 3$.

In \cite{Hoff} Hoffman proved that on the infinite cluster in $\Z^2$
there are flows with finite $(2,2+\ep)$-energy. We prove the following:

\begin{thm} \label{thm:z2_flow}
The infinite cluster of super critical bond percolation in $\Z^2$
a.s.\ supports
a flow of finite $(2,1+\ep)$-energy.
\end{thm}

Theorem \ref{thm:z2_flow} is a corollary of a more general result.

\begin{thm}\label{thm:hoffgen}
Let $g:[0,\infty]\to[0,\infty]$ be a convex function s.t.
\\(I) There exits $l\in\N$ such that $x^{-l}g(x)$ is decreasing.
\\(II) $\Z^2$ supports a flow with finite $g$-energy.
\\Then, for every $p>\frac{1}{2}$, a.s. the infinite cluster for
bond percolation (with parameter $p$) in $\Z^2$,
 supports a flow of finite $g$-energy.
\end{thm}

In section \ref{sec:firstproof} we give a proof of theorem \ref{thm:core}
and show that it implies theorem \ref{thm:main}. The argument is based on
large deviation results by Antal-Pisztora (\cite{AnPi}).

In section \ref{sec:secondproof} we give an alternative proof of Theorem
\ref{thm:main} for large values of the percolation parameter $p$. This
proof relies on connectivity properties of $\Z^d$ instead of the
Antal-Pisztora theorem and can be extended to other graphs.
This alternative approach yields a new proof of
the following theorem due to Benjamini and Schramm \cite{BeSc}:

\begin{thm} \label{thm:cayley}
Let $G$ be the Cayley graph of a finitely generated group of
polynomial growth which is not a finite extension of $\Z$ or
$\Z^2$. Then for $p$ sufficiently close to $1$, the infinite
$p$-percolation cluster on $G$ is transient.
\end{thm}

In section \ref{sec:z2_flow} we
prove Theorems \ref{thm:hoffgen} and \ref{thm:z2_flow}.

\section{Proof of Theorem \ref{thm:core}} \label{sec:firstproof}

We need the following lemma (See e.g. chapter $2$ of \cite{yvkl}
or \cite{climb}):
\begin{lemma} \label{lem:yvkl}
Let $G$ be a graph, and let $v_0\in G$. Then, $G$ is transient if and only
if there exists a probability measure $\mu$ on self-avoiding paths in $G$
starting at $v_0$ s.t. if $P$ and $Q$ are two paths chosen independently
according to $\mu$, the expected number of edges in $P\cap Q$ is finite.
\end{lemma}

We identify $\mu$ with the flow that it induces. Note that the condition on
$\mu$ is equivalent to the condition:
\begin{equation} \label{eq:finiteinter}
H_{2,0}(\mu) = \sum_e \mu(\{P | e \in P\})^2 < \infty
\end{equation}

We also need the next lemma which is what is actually proved (even though
not stated) by Antal and Pisztora in \cite{AnPi}. It can be found in this
form as Lemma 2.13 of \cite{GLA}.

{}For two vertices in the same connected component, $x$ and $y$, denote by
$D(x,y)$ the length of the shortest path between them.

\begin{lemma}
Let $p>p_c(\Z^d)$. There exist $\rho = \rho(p,d)$ and $C>0$ such that for
every integer $m$ and any $v, w \in \Lambda_m:=[-m,m]^d$:
\begin{equation} \label{eq:3.27}
\P_p(v \leftrightarrow w \ {\rm and} \ D(v,w)>\rho m) \leq e^{-C m}
\end{equation}
\end{lemma}

We will use the following equivalent form of this lemma:
\begin{lemma} \label{lem:gla}
If $p>p_c(\Z^d)$ then there exist $\rho = \rho(p,d)$ and $\theta>0$ such
that if $v,w$ are both in the infinite cluster then for all $m>\rho d(v,w)$
\begin{equation*}
\P_p( D(v,w) > m ) \leq e^{-\theta m}
\end{equation*}
\end{lemma}

The idea of the proof of Theorem \ref{thm:core} is to construct a measure
on paths in the infinite percolation cluster by taking the measure $\mu$
supported on paths in the core, and modifying the paths to be in the
percolation cluster. A path chosen by $\mu$ will typically have gaps where
it passes through closed edges or non-percolating vertices. We will replace
the gap in the path by a ``bridge'' which will be the shortest path in the
infinite cluster which connects the ends of the gap. We will see that with
positive probability all the gaps that are inside the $C$-core of $G$ will
have bridges in $G$.

\begin{defn}
A {\bf gap} in a percolation configuration is a connected cluster in the
complement of the infinite cluster.
\end{defn}

Clearly a gap consists of a closed cluster and anything separated by it
from the infinite cluster. Therefore if in some configuration large closed
clusters are rare then so are large gaps. This will be stated formally later.

We first proceed to prove Theorem \ref{thm:core} for $p$ close to $1$:

\begin{lemma}\label{lem:high_p}
Let $d\geq3$. Then there exist  $p_d<1$ and $C=C_d$ with the following
property: Let $p\geq p_d$, and let $A$ be a subgraph $\Z^d$ s.t.
  (1) the $C$-core of $A$ is transient, and
  (2) $p$-percolation on $A$ has, a.s., a unique infinite cluster.
Then, a.s., the infinite $p$-percolation cluster of $A$ is transient.
\end{lemma}

\begin{proof}
Choose $p_d$ close enough to $1$, so that $(1-p_d)$-percolation is
sub-critical. Consider percolation on $A$ as a restriction to $A$ of
Bernoulli percolation on all of $\Z^d$. Denote by $I$ the infinite
percolation cluster in $A$, and by $J \supseteq I $ the infinite
percolation cluster in $\Z^d$.

A.s., $A\setminus J$ has only finite clusters. Moreover, if for a vertex
$x$ we use $C(x)$ to denote the cluster of $A\setminus J$ containing $x$,
then there exist $\delta$ and $\alpha>0$ such that for every $x$,
\begin{equation}\label{eq:pregamma}
\P(|C(x)|>n)<\delta e^{-n^\alpha}.
\end{equation}
{}For the diameter we even have an exponential bound: There exists a
constant $\gamma$ such that
\begin{equation}\label{eq:gamma}
\P(\diam(C(x))>n)<e^{-\gamma n}.
\end{equation}

Choose
\begin{equation*}
C_d=l=2\max\left(\frac{d\rho}{\gamma},\frac{d\rho}{\theta}\right)
\end{equation*}
where $\rho$ and $\theta$ are from Lemma \ref{lem:gla} and $\gamma$ is from
(\ref{eq:gamma}). The $l$-core of $A$ is transient. Lemma \ref{lem:yvkl}
states that there is a probability measure $\mu$, satisfying
(\ref{eq:finiteinter}), on paths in the core starting at some $v_0$. Since
transience of $I$ is a 0-1 event, and $v_0 \in I$ with positive
probability, we may assume that $v_0\in I$.

Let $P=(v_0,v_1,v_2,v_3,...)$ be a path (chosen according to $\mu$).
A.s., $P$ intersects $J$ infinitely often, so we may restrict ourselves to
paths that intersect $J$ infinitely often. We modify $P$ to get $P'$, a
path in $J$, as follows: At any time at which $P$ enters a gap $B$, we will
replace the part of the path which is in $B$ by a shortest possible path,
inside the cluster, between its ends. Let $\mu'$ be the new measure on
paths (or, more precisely, if $\phi$ is the function which assigns $P'$ to
each $P$, then let $\mu'=\mu\circ\phi^{-1}$).

Clearly $\mu'$ is supported on paths in $J$ starting at $v_0$. To conclude
the proof of the transience of $I$ we will prove the following two lemmas:

\begin{lemma} \label{lem:mu_support}
In the above setting, for a suitable value of $l$, with positive
probability $\mu'$ is supported on paths in $I$.
\end{lemma}

\begin{lemma} \label{lem:mu_energy}
In the above setting, a.s. $\mu'$ has finite energy (i.e. satisfies
{\it {\em (\ref{eq:finiteinter})}}).
\end{lemma}

{}From these two lemmas it follows that with positive probability we have
constructed a measure on paths in $I$ satisfying (\ref{eq:finiteinter}), so
with positive probability $I$ is transient. Since transience of $I$ is a
0-1 event we are done.
\end{proof}   

%
\begin{proof} [Proof of Lemma \ref{lem:mu_support}]
We will show that with positive probability the paths of $\mu'$ are in $A$
and hence in $I$.

Let $\G_R$ be the set of gaps at distance $R$ from the origin. Trivially
$|\G_R| < C_1 R^{d-1}$.

Let $A_R(\beta), \beta>0$ be the event that all gaps in $\G_R$ have
diameter less then $\beta d \log R$. By the bound on the distribution of
the diameter of the gaps,
\begin{equation*}
\P(A_R(\beta)) \geq 1 - (C_1 R^{d-1}) e^{-\gamma \beta d \log R}
         =  1 - \frac{C_2}{R^{1+d(\gamma \beta-1)}}
\end{equation*}

Choose $\beta > \max\left(\frac{1}{\gamma},\frac{1}{\theta}\right)$.
Since $A_R(\beta)$ are monotone increasing events and since
$\gamma \beta > 1$ we have
\begin{equation*}
\P(\bigcap_R A_R(\beta))  \geq  \prod_R \P(A_R(\beta))  >  0
\end{equation*}

{}From now on we assume that $\bigcap_R A_R(\beta)$ occurs. A gap of
diameter $s$ has less then $s^d$ vertices in it, so the number of bridges
over such a gap is less then ${s^d} \choose 2$. Assuming $A_R$, the total
number of bridges over gaps in $\G_R$ is bounded by
$(C_1 R^{d-1}) (\beta d \log R)^{2d}$.

Consider a bridge $b$ over a gap in $\G_R$. Its endpoints are no further
then the diameter of the gap, or $\beta d \log(R)$. By Lemma \ref{lem:gla},
since $l$ is chosen s.t. $l / (\beta d) > \rho$ then
\begin{eqnarray*}
\P(\len(b) > l \log R)
  &<& \P(\len(b) > \rho \beta d \log R)              \\
  &<& e^{-\theta \beta d \log R }                    \\
  &=& R^{-\theta \beta d}
\end{eqnarray*}

The probability of the event $\Gamma_R$ that all of the bridges over gaps
in $\G_R$ are shorter then $l \log R$ is at least
\begin{equation*}
1- (C_1 R^{d-1}) (\beta d \log R)^{2d} R^{- \theta \beta d}
 = 1 - C_3 \frac{(\log R)^{2d}}{R^{1 + d(\theta \beta - 1)}}
\end{equation*}

Since $\{\Gamma_R\}_{R\geq 1}$ and $\{A_R\}_{R\geq 1}$ are
increasing events, conditioning on $\bigcap_R A_R$ can only
increase the probability of $\bigcap_R \Gamma_R$. Since $\beta$ is
s.t. $\theta \beta > 1$ and $\gamma \beta > 1$, we get that with
positive probability this occurs for every $R$.

A bridge over a gap in $\G_R$ in the $l$-core of $A$, of length at most
$l\log R$ cannot leave $A$, and so we are done.
\end{proof} 


%

\begin{proof} [Proof of Lemma \ref{lem:mu_energy}]
Define the functions
\begin{equation*}
{}F(e) = \mu(\text{paths which go through e})
\end{equation*}
and
\begin{equation*}
{}F'(e) = \mu'(\text{paths which go through e})
\end{equation*}
The $\mu$-expected (resp. $\mu'$-expected) size of the intersection between
independently chosen paths is exactly $\sum_e{F(e)^2}$
(resp. $\sum_e{F'(e)^2}$).

We now show that
\begin{equation*}
\E{\left(\sum_{e}{F'(e)^2}\right)}<\infty
\end{equation*}
and thus a.s. $\sum_e{F'(e)^2} < \infty$.

We say that an edge $e$ is projected on an edge $f$ if either $e=f$ or if
some path with a gap including $e$ has a bridge over it passing through $f$.
We will denote the event that $e$ is projected on $f$ by $e \to f$. It is
clear that
\begin{equation*}
{}F'(f)\leq\sum_{e \to f}{F(e)}
\end{equation*}

Define $S(e) = \{f | e \to f \}$ and $T(f) = \{e | e \to f\}$.
Then we have the bound:
\begin{eqnarray*}
\sum_{f\in E} F'(f)^2
  &\leq& \sum_{f\in E} \left( \sum_{e\in T(f)} F(e) \right)^2            \\
  &\leq& \sum_{f\in E} |T(f)| \sum_{e\in T(f)} F(e)^2                    \\
  & =  & \sum_{e\in E}\left(F(e)^2\cdot\sum_{f\in S(e)}|T(f)|\right)
\end{eqnarray*}

We use the following claim:
\begin{claim}\label{numsol}
There exist constants $C$ and $a>0$ s.t. for every $e$ and $f$,
\begin{equation}\label{eq:numsol1}
\P(e\to f) < C\exp(-a d(e,f)).
\end{equation}
\end{claim}


By (\ref{eq:numsol1}) we see that there exist $\alpha$ and $\beta$ s.t.
\begin{equation}\label{eq:numsol2}
\P(|S(e)| > n) < e^{-\alpha n^\beta}
\end{equation}
and
\begin{equation}\label{eq:numsol3}
\P(|T(f)| > n) < e^{-\alpha n^\beta}.
\end{equation}

{}For every $e$ and $g$,
\begin{eqnarray*}
\P(\exists_f(f\in S(e) \ve g\in T(f)))
   &\leq&  \sum_f\P(f\in S(e) \ve g\in T(f))            \\
   &\leq&  e^{-\alpha'' d(e,g)^{\beta''}}.
\end{eqnarray*}

for appropriate $\alpha''$ and $\beta''$. Therefore,
\begin{eqnarray*}
\P\left(\#\left(\bigcup_{f\in S(e)}T(f)\right)>n\right)
< e^{-\alpha' n^{\beta'}}
\end{eqnarray*}
for some $\alpha'$ and $\beta'$. Since
\begin{eqnarray*}
\sum_{f\in S(e)}|T(f)|  \leq
\#\left(\bigcup_{f\in S(e)}T(f)\right)  \cdot  |S(e)|
\end{eqnarray*}
we get that
\begin{eqnarray*}
\P(\sum_{f\in S(e)}|T(f)|>n)
  &\leq& \P(|S(e)|>n^{1/2}) +
         \P\left(\#\left(\bigcup_{f\in S(e)}T(f)\right)>n^{1/2}\right)  \\
  &\leq& e^{-\alpha n^{\beta/2}} + e^{-\alpha' n^\frac{\beta'}{2}}.
\end{eqnarray*}


In particular, for some constant $D<\infty$
\begin{equation}\label{eq:pre_few_sols}
\E\left(\sum_{f\in S(e)}|T(f)|\right)<D
\end{equation}
and therefore
\begin{equation*}
\E \left(\sum_f F'(f)^2\right) < \infty
\end{equation*}
as required.
\end{proof}    

\begin{remark}\label{eq:few_sols}
Actually, the same estimate implies that every moment of
\begin{equation*}
\sum_{f\in S(e)}|T(f)|
\end{equation*}
is bounded.
\end{remark}

\begin{proof}[Proof of Claim \ref{numsol}]
Let $n=d(e,f)$. We assume that $e\not\in I$ (because otherwise $e \to f$
only if $e=f$). Let $\Upsilon$ be the connected component of $\Z^d-J$
containing $e$.
\begin{equation}\label{upsilon}
\P(e\to f) \leq \P\left(|\Upsilon|>\sqrt{n}\right) +
\P\left(e\to f \big| |\Upsilon|\leq\sqrt{n}\right).
\end{equation}
Since $\Z^d-J$ is the union of sub-critical percolation and finite clusters
of supercritical percolation, for suitable $C$ and $\alpha$,
\begin{equation}\label{upsilon1}
P \left(|\Upsilon|>\sqrt{n}\right) \leq Ce^{-n^\alpha}
\end{equation}
To bound the second term in the RHS of (\ref{upsilon}), note that
$\P\left(e\to f \big | |\Upsilon|\leq\sqrt{n} \right)$ is bounded by the
probability that there exist $a$ and $b$, each at distance at most
$2\sqrt{n}$ from $e$, s.t. the Antal-Pisztora path between $a$ and $b$ goes
through $f$. This is bounded by
\begin{eqnarray}\label{upsilon2}
\nonumber
\sum_{a,b}
\P\left(f\in\text{bridge}(a,b)\big||\Upsilon|\leq\sqrt{n}\right)
&\leq&
\sum_{a,b}\P\left(D(a,b)\geq n\big||\Upsilon|\leq\sqrt{n}\right)   \\
\leq
\sum_{a,b}\P\left(D(a,b)\geq n\right)
&\leq&
4^d n^{d/2} e^{-\frac{\theta}n}
\end{eqnarray}
for $n$ large enough, where the second inequality is due to FKG, and the
third one is due to Lemma \ref{lem:gla}. (\ref{upsilon}), (\ref{upsilon1})
and (\ref{upsilon2}) now give the desires result.
\end{proof}

In order to move to lower values of $p$, we need to use a renormalization
technique. We use the same technique used by \olle and Mossel in \cite{EO}:
Notice that Lemma \ref{lem:high_p} is also valid when we consider site
percolation with high enough retention probability. Let $\hat{p}<1$ be such
a high enough retention probability. Let $p\in(p_c,1]$.  Let $N$ be a
(large) positive integer, divisible by $8$. Let $A$ be a subgraph of $\Z^d$
such that its $CN$-core is transient for $C=C_d$ and $N$ as defined. Define
\begin{equation*}
A_N=N^{-1}\left\{v\in N \cdot \Z^d \big|
v+\left[-\frac{5N}{8},\frac{5N}{8}\right]^d\subseteq A\right\}.
\end{equation*}
$A_N$ has a transient $C$-core, and therefore there is (a.s.) a transient
infinite cluster in the $\hat{p}$-site percolation on $A_n$. For each $v$
in $N\cdot A_N$ define $Q_N(v)$ to be the cube of side-length $5N/4$
centered at $v$. Apply $p$-percolation on $A$, and let $A_p$ be the random
set of vertices $v\in N\cdot A_N$ such that $Q_N(v)$ contains a connected
component which connects all $2d$ faces of $Q_N(v)$ but contains no other
connected component of diameter greater than $N/10$. It follows from
Proposition 2.1 in Antal-Pisztora \cite{AnPi} that if $N$ is large enough
then $A_p$ dominates the $\hat{p}$-site percolation on $A_N$. Therefore the
$p$-percolation on $A$ has a transient infinite cluster.

In order to prove Theorem \ref{thm:main} we only have to prove the
following lemma:
\begin{lemma}\label{lem:yeshcore}
Let $h$ be monotone s.t. $\W_h$ is transient. Then, for every $C$, the
$C$-core of $\W_h$ (Denoted by $\W(h,C)$) is transient as well.
\end{lemma}

\begin{proof}
Recall Lyons' criterion for transience of a wedge: $\W_h$ is transient if
and only if
\begin{equation*}
\sum_{j=1}^\infty  \frac{1}{jh(j)}  <  \infty
\end{equation*}
Lyons also proved that if the wedge $\W_h$ is transient, then its subgraph
\begin{equation*}
\V_h=\W_h\cap\{(x,y,z)|x\geq 0\text{ and }|y|\leq x\}
\end{equation*}
is also transient.

Now, let $C$ be arbitrary and let $h$ be an increasing function such that
$W_h$ is transient.

\begin{claim}\label{wlog}
We may assume, W.l.o.g, that $h(x+1)-h(x) \leq 1$ for every $x$.
\end{claim}

Lyons' condition (\ref{eq:lyons}) implies that There exists $x_0$ s.t. for
all $x>x_0$, we have $h(x) > 4C\log(x)$. Define
\begin{equation*}
g(x)=\frac{1}{2}h(x + x_0).
\end{equation*}
We will show that
  \\(A) $W(h,C)$ contains a translation of $\V_g$, and that
  \\(B) $g$ satisfies Lyons' condition (equation (\ref{eq:lyons})).
\\This will suffice for proving the lemma.

{}For (A), we just need to prove that $\V_g+(x_0,0,0)\subseteq\W(h,C)$.
Let $a=(x,y,z) \in \V_g+(x_0,0,0)$. The norm of $a$ is
\begin{equation*}
(x^2+y^2+z^2)^\frac{1}{2}\leq 2x
\end{equation*}
Therefore, we have to show that every vertex at distance up to $C\log(2x)$
from $a$ is in $W_h$, but this is clear because if $a_1=(x_1,y_1,z_1)$
satisfies $|x-x_1|\leq C\log(2x)$ and $|z-z_1|\leq C\log(2x)$ then
\begin{eqnarray*}
z_1 \leq g(x)+C\log(2x)
  & = & g(x)+C\log(x)+C\log(2)                          \\
  &\leq& g(x)+C\log(x_1)+2C\log(2)                      \\
  &\leq& g(x_1) + 2C\log(x_1) + 2C\log(2)               \\
  &\leq& h(x_1)
\end{eqnarray*}

(B) is trivial.
\end{proof}

\begin{proof}[Proof of Claim \ref{wlog}]
Given $h$, define $f$ inductively to be: $f(0)=h(0)$ and
\[
f(n+1)=\min(h(n+1),f(n)+1).
\]
It suffices to show that if
\begin{eqnarray*}
U(h)=\sum_{j=1}^\infty  \frac{1}{jh(j)}  <  \infty
\end{eqnarray*}
then
\begin{eqnarray*}
U(f)=\sum_{j=1}^\infty  \frac{1}{jf(j)}  <  \infty.
\end{eqnarray*}
Let $A=\{j:f(j)=h(j)\}$ and $B=\{j|f(j)=f(j-1)\}$. Then, $B \subseteq A$.
Now,
\begin{equation}\label{wlog1}
\sum_{j\in A} \frac{1}{jf(j)} = \sum_{j\in A} \frac{1}{jh(j)} \leq
U(h) < \infty,
\end{equation}
and
\begin{equation}\label{wlog2}
\sum_{j\not\in A} \frac{1}{jf(j)}\leq\sum_{j\not\in B}
\frac{1}{jf(j)}\leq  \sum_{j=1}^\infty j^{-2} <\infty.
\end{equation}
Combining (\ref{wlog1}) and (\ref{wlog2}) we get that $U(f)<\infty$.
\end{proof}

\section{General Graphs} \label{sec:secondproof}

If the percolation parameter $p$ is close enough to $1$,
specifically if $p>1-p_c$, then a.s. the closed edges compose only
finite closed clusters. If $p>p_s$ for some (larger) critical
$p_s$ then the closed clusters will not only be finite, but they
will also be isolated from each other, and then a path that passes
through a gap $C$ can be bridged without leaving its boundary
$\partial C$. This notion will lead us to another proof for
``nice'' graphs.

\begin{defn}
A graph $G$ is said to have {\bf $k$-connected boundaries} if for any connected
set of vertices $A$ s.t. $G \setminus A$ is also connected, the subgraph
spanned by
\[
B_k(A) = \left\{ v | D(v,A) < k \right\} \setminus A
\]
is connected. (i.e. the set of vertices outside $A$ but at distance at
most $k$ from $A$ is connected).
\end{defn}

The property of having $k$-connected boundaries is rather general and
holds for many graphs. For example any planar graph $G$ which can be
embedded in the plane with all faces having at most $2k+1$ edges, has
$k$-connected boundaries, since if a set $A$ includes a vertex from a
face then all other vertices of that face will be within distance $k$
from $A$. (See for example Lemma 4.4 of \cite{HeSc}).

A special case of a result in \cite{BaBe} tells us that minimal cut
sets in $\Z^d$ are 2-connected (i.e. together with their neighbors
form connected sets). It follows from this that $\Z^d$ has 3-connected
boundaries.




We say that an edge is k-{\bf strongly open} if it is open and so are
all the edges at distance up to $k$ from it. An edge is
k-{\bf weakly closed} if it is not k-strongly open.

\begin{lemma} \label{lem:seperation}
If a graph $G$ has degrees bounded by $d$, there exists
$p_s=p_s(k)$ such that for $p>p_s$ a.s. the k-weakly closed edges of
Bernoulli percolation configuration do not percolate and s.t. if
$C(v)$ is the weakly closed cluster of vertex $v$ then for
appropriate $\gamma=\gamma(p)$
\[
\P(\diam(C(v)) > n)  <  e^{-\gamma n}
\]
\end{lemma}

\begin{proof}
Let $d$ be a bound on the degrees of $G$. Since whether an edge is
k-weakly closed depends only on those edges up to distance $k$ from it,
the configuration of k-strongly open edges
dominates bond-percolation with parameter $q=q(p)$, and moreover
$q(p) \to 1$ as $p \to 1$ (See \cite{dompq}). If $p>p_s$ such that
$q(p_s) > 1 - d^{-2}/5$ then the weakly closed clusters are dominated
by sub-critical Galton-Watson trees. Therefore there exists a
$\gamma>0$ and $w>5$ s.t.
\begin{equation}\label{bddiam}
\P(\diam(C(v)) > n)  <  e^{-\gamma n}
\end{equation}
and
\begin{equation}\label{bdvol}
\P(|C(v)| > n)  <  c \cdot n^{-w}
\end{equation}
\end{proof}

\begin{defn}
The {\bf $k$-inner boundary} of a set $A$ in a graph $G$ is the
$k$-boundary of $G-A$, i.e.
$
\{x\in A: d(x,G \setminus A)<k\}
$
\end{defn}

\begin{lemma}\label{innbnd}
Let $G$ be a graph with $k$-connected boundaries. Let $\bar{B}$ be a
connected component of the complement of the k-strongly open cluster. Let $U$
be the $k$-inner boundary of $\bar{B}$. Then
\\(A) $U$ is in the open cluster.
\\(B) $U$ is connected.
\end{lemma}

\begin{proof}
(A) follows from the definition of the inner boundary and of k-strongly
open edges. (B) follows from the fact that $U$ is the $k$-boundary of
$G-\bar{B}$ and the fact that $G$ has $k$-connected boundaries.
\end{proof}

\begin{thm} \label{thm:good_graph}
Let $G$ have bounded degrees and $k$-connected boundaries. Let $A$
be a sub-graph of $G$ such that for every $C$, the $C$-core of $A$
is transient. Assume also that for $p$ close enough to $1$, the
$p$-percolation on $A$ has, a.s., a unique infinite cluster. Under
those conditions, if $p$ is close enough to $1$, then the infinite
$p$-percolation cluster of $A$ is transient.
\end{thm}

\begin{proof}[Proof of Theorem \ref{thm:good_graph}]
{}For $p>p_s$ where $p_s$ is as defined in Lemma \ref{lem:seperation}.
Consider percolation on $A$ as a restriction to $A$ of Bernoulli
percolation on all of $G$. Denote by $I$ the infinite percolation cluster
in $A$, and by $J \supseteq I$ the infinite percolation cluster in $G$.
If $C(v)$ is the k-weakly closed cluster containing $v$ then we know that
$\P(\diam(C(v)) > n) < e^{-\gamma n}$.

As before, For $l$ (to be determined later) the $l$-core of $A$ is
transient. Lemma \ref{lem:yvkl} states that there is a probability measure
$\mu$, satisfying (\ref{eq:finiteinter}), on paths in the core starting at
some $v_0$. Since transience of $I$ is a 0-1 event, and $v_0 \in I$ with
positive probability, we may assume that $v_0\in I$.

Let $P=(v_0,v_1,v_2,v_3,...)$ be a path (chosen according to $\mu$).
A.s., $P$ intersects $I$ infinitely many times, so we may restrict
ourselves to such paths. Now, we modify $P$ to get $P'$, a path in $J$, as
follows: At any time at which $P$ enters a gap $B$, we consider the k-weakly
closed extension of the gap, $\bar B$, and  we replace a part of the path from
the first time $p$ reaches $\bar B$ until the last time it leaves it, by a
path in the inner $k$-boundary of $\bar B$.

By Lemma \ref{innbnd}~(B), such a path exists. If $\phi$ is the function
which assigns $P'$ to each $P$, then let $\mu'=\mu\circ\phi^{-1}$. By part
(A) of the same lemma, $\mu'$ is supported on paths in $J$.

The length of a bridge over a gap of size $n$ is bounded by the size of the
gap's $k$-boundary.
Using the estimates (\ref{bddiam}) and (\ref{bdvol}), we can repeat the
proofs of Lemma \ref{lem:mu_support} and Lemma \ref{lem:mu_energy} and get
that with a positive probability, $\mu'$ is a measure on paths in $I$ with
finite energy, which proves the transience of $I$.
\end{proof}

In the case where $A$ is the whole graph $G$ this takes the form:

\begin{cor}
Under the conditions of Theorem \ref{thm:good_graph}, if $G$ is transient
then the $p$-percolation cluster in $G$ is transient.
\end{cor}

Now we can prove Theorem \ref{thm:cayley}.
\begin{proof}[Proof of Theorem \ref{thm:cayley}]
By Gromov's theorem, (see \cite{GRO}) any finitely generated group
of polynomial growth is a finite extension of a finitely presented
group. Corollary 4 of \cite{BaBe} tells us that the Cayley graph
of such a group has $k$-connected boundaries for some $k$
(depending on the lengths of the relations). Corollary 10 of the
same paper tells us that for $p$ close enough to $1$, the infinite
cluster is unique. Since $G$ is not a finite extension of $\Z$ or
$\Z^2$, by a theorem of Varopoulos, $G$ is transient (See
\cite{varop}). Thus all the conditions of Theorem
\ref{thm:good_graph} are satisfied.
\end{proof}

\section{Proof of Theorem \ref{thm:z2_flow}} \label{sec:z2_flow}

We first restate and prove Theorem \ref{thm:hoffgen}.
\\{\bf Theorem \ref{thm:hoffgen}.}
{\em Let $\varphi:[0,\infty]\to[0,\infty]$ be a convex function
s.t.
\\(I) There exits $l\in\N$ such that $x^{-l}\varphi(x)$ is decreasing.
\\(II) $\Z^2$ supports a flow with finite $\varphi$-energy.
\\Then, for every $p>\frac{1}{2}$, a.s. the infinite percolation cluster
supports a flow with finite $\varphi$-energy.   }

\begin{proof}
Throughout the proof the term energy will refer to
$\varphi$-energy. We begin with a flow $F$ on $\Z^2$ with finite
energy with the source at 0. Such a flow exists by condition (II)
of the theorem. The function $\varphi$ is increasing and convex in
the flow on each edge, so the flow $F$ can be made acyclic (i.e. a
flow s.t. there is no cycle s.t. all of its edges get positive
flow) without increasing its energy. Normalize the flow so that
the total flow out of 0 is 1 and  the energy is finite. Since $F$ is
an acyclic flow with source at 0, it induces a probability measure
$\mu$ on
self-avoiding paths starting at 0, see e.g. \cite{resist}.
The $\mu$-measure of the set of paths passing
through $e$ is exactly $|F(e)|$.

As before, note that the existence of a flow with finite energy on the
percolation cluster is a 0-1 event, so we can assume that 0 belongs to the
infinite cluster. Since in $\Z^2$ we have $p_c=\frac{1}{2}$ our condition
$p>1-p_c$ is equivalent to $p>p_c$. We construct as before a new measure
$\mu'$, by using the shortest bridge over any gap in a path. We wish to
show that $\mu'$ a.s. has finite energy.

The function $\varphi$ is increasing and convex in $[0,1]$, and,
moreover,
\begin{eqnarray} \label{eq:h_bound}
\varphi \left( \sum_{i=1}^n x_i \right) & = & \varphi\left( n^{-1}
\sum_{i=1}^n nx_i \right)
\nonumber \\
& \leq & n^{-1}\sum_{i=1}^n\varphi(n x_i)
\nonumber \\
& \leq & n^{-1}\sum_{i=1}^nn^l\varphi(x_i)
\nonumber \\
& = & n^{l-1} \sum_{i=1}^n \varphi(x_i)
\end{eqnarray}

We now proceed to estimate the energy of $F'$. Using the notation in the
proof of Lemma \ref{lem:mu_energy}:
\begin{eqnarray*}
\varphi(F'(f)) & \leq &  \varphi(\sum_{e\to f} F(e))              \\
         & \leq &  |S(f)|^{l-1} \sum_{e\to f} \varphi(F(e))
\end{eqnarray*}

Summing over $f$ we get:
\begin{eqnarray*}
H_\varphi(F)
  & = &  \sum_f\varphi(F'(f))                                 \\
  &\leq& \sum_f |S(f)|^{l-1} \sum_{e\to f} \varphi(F(e))      \\
  & = &  \sum_e \varphi(F(e)) \times
            \left(\sum_{f\in S(e)}|T(f)|^{l-1}\right)   \\
  &\leq& \sum_e \varphi(F(e)) \times
            \left(\sum_{f\in S(e)}|T(f)|\right)^{l-1}
\end{eqnarray*}

By Remark \ref{eq:few_sols} we have
\[
\E H_\varphi(F') \leq \E \left[
         \left(\sum_{f\in S(e)}|T(f)|\right)^{l-1}
\right]
 H_\varphi(F) < \infty
\]

And thus $F'$ a.s. has finite energy.
\end{proof}

In order to prove Theorem \ref{thm:z2_flow}, we just have to notice that
the function
$\Psi_{2,1+\ep}=\frac {|x|^2} {\log (1+|x|^{-1})^{1+\ep}}$ corresponding
to $(2,1+\ep)$-energy satisfies the conditions of Theorem \ref{thm:hoffgen}
with $l=4$.


\noindent
Omer Angel and Itai Benjamini,\\
The Weizmann Institute and Microsoft research\\
omer@wisdom.weizmann.ac.il \\
itai@wisdom.weizmann.ac.il\\
\\
Noam Berger and Yuval Peres,\\
The University of California\\
noam@stat.berkeley.edu \qquad
\\
peres@stat.berkeley.edu\\
\end{document}